\documentclass[twocolumn,10pt]{asme2ej}
\usepackage{epsfig} 
\usepackage{amsmath,amssymb,,amsbsy,amsfonts}
\usepackage{graphics}
\usepackage{float}
\usepackage{caption}
\usepackage{graphicx}
\usepackage{multirow}
\usepackage{mathtools}
\usepackage{subfig}
\usepackage{tikz}
\usepackage{nomencl}
\usepackage{cleveref}
\usepackage{cite}
\usepackage{epsfig} 
\usepackage{fancyhdr}
\usepackage{color}
\usepackage{float}

\usepackage{url}
\usepackage{breqn}
\usepackage[toc,page,header]{appendix}
\usepackage{multicol}

\title{Adaptive Discrete Second Order Sliding Mode Control with Application to \\Nonlinear Automotive Systems}

\author{Mohammad Reza Amini\thanks{Address all correspondence to this author.}\\
 {\tensfb Mahdi Shahbakhti}
   \affiliation{Mechanical Eng.-Eng. Mechanics Dept.\\
    Michigan Technological University\\
    Houghton, Michigan 49931\\
    Emails: mamini, mahdish@mtu.edu
    }	
}

\author{Selina Pan
    \affiliation{Department of Mechanical Engineering\\
    University of California\\
    Berkeley, CA 94720\\
    Email: slpan@berkeley.edu
    }
}

\author{\vspace{-2cm}}
\graphicspath{{Figures/}}
\begin{document}
\maketitle
\begin{abstract}
\noindent\textbf{Abstract}\\
{\it Sliding mode control (SMC) is a robust and computationally efficient model-based controller design technique for highly nonlinear systems, in the presence of model and external uncertainties. However, the implementation of the conventional continuous-time SMC on digital computers is limited, due to the imprecisions caused by data sampling and quantization, and the chattering phenomena, which results in high frequency oscillations. One effective solution to minimize the effects of data sampling and quantization imprecisions is the use of higher order sliding modes. To this end, in this paper, a new formulation of an adaptive second order discrete sliding mode control (DSMC) is presented for a general class of multi-input multi-output (MIMO) uncertain nonlinear systems. Based on a Lyapunov stability argument and by invoking the new Invariance Principle, not only the asymptotic stability of the controller is guaranteed, but also the adaptation law is derived to remove the uncertainties within the nonlinear plant dynamics. The proposed adaptive tracking controller is designed and tested in real-time for a highly nonlinear control problem in spark ignition combustion engine during transient operating conditions. The simulation and real-time processor-in-the-loop (PIL) test results show that the second order single-input single-output (SISO) DSMC can improve the tracking performances up to 90\%, compared to a first order SISO DSMC under sampling and quantization imprecisions, in the presence of modeling uncertainties. Moreover, it is observed that by converting the engine SISO controllers to a MIMO structure, the overall controller performance can be enhanced by 25\%, compared to the SISO second order DSMC, because of the dynamics coupling consideration within the MIMO DSMC formulation.}
\end{abstract}
\vspace{-0.6cm}

\section{Introduction}
\label{sec:intro}
\label{Sec:Intro}
Converting a high dimensional tracking control problem into a low dimensional stabilization control problem is the key feature of sliding mode control (SMC)~\cite{Slotine,Amini_ACC2018}. SMC shows robust characteristics against external disturbances and model uncertainty/mismatch, while requiring low computational efforts. However, there are challenging issues that arise during implementation of SMC on digital processors, which limit the real-time application of SMC. The two well-recognized challenging issues include: (i) high frequency oscillations due to chattering phenomenon, and (ii) implementation imprecisions due to the analog-to-digital (ADC) converter unit~\cite{utkin2013sliding,Amini_DSC,AminiSAE2016}. 

The concept of higher order sliding modes for continuous-time systems is shown to be an effective approach for reducing the oscillation due to chattering. This approach was first introduced in the 1980s~\cite{nollet2008observer}. The basic idea of the higher order SMC is to not only steer the sliding function to the sliding manifold, but also drive all the higher order derivatives of the sliding variable to zero. Higher order SMC reduces the high frequency oscillations by transferring the chattering caused by the discontinuity to the higher order sliding mode derivatives~\cite{nollet2008observer}. Higher order SMC leads to less oscillations; however, it adds complexity to the calculations. Moreover, it has been shown in~\cite{Acary,huber2016implicit} that converting the continuous-time SMC to a discrete sliding mode controller (DSMC), by using an implicit Euler discretization, allows for a drastic decrease in the chattering in both the input and the output. Thus, according to~\cite{mihoub2009real}, which presents a second order DSMC, the idea of higher order DSMC can be an ideal solution for the chattering problem by taking advantage of characteristics of the higher order SMC and discretized SMC. 

In addition to the high frequency oscillations issue, a gap often occurs between the designed and the implemented conventional SMCs, which degrades the controller performance from its expected behaviour significantly~\cite{Amini_CEP,Amini_DSC}. This gap is mostly created due to (i) data sampling and quantization imprecisions that are introduced by the ADC at the controller input/output (I/O), and (ii) uncertainties in the modeled dynamics. Fortunately, the SMC structure allows for further modification to improve the controller robustness against ADC (sampling and quantization) imprecisions, and compensate for the uncertainties within the plant model~\cite{Amini_CEP}. 

{There are several works in the literature aimed at improving the robustness of conventional SMC design against ADC and modeling uncertainties.} The results for SMC with incorporated maximum ADC uncertainty bounds~\cite{KyleACC}, and uncertainty adaptive SMC~\cite{AminiSAE2016} have shown that the continuous-time SMC performance can be improved against implementation imprecisions. However, these improvements are limited, and SMC fails at higher sampling times. On the other side, it has been shown in~\cite{Kyle_DSCC} that the DSMC and DSMC with incorporated maximum ADC uncertainty bounds can significantly improve the controller robustness against sampling and quantization imprecisions. In the recent works~\cite{Amini_ACC2016,Amini_DSCC2016}, we have shown that by incorporating an online ADC uncertainty prediction and propagation mechanism, not only the robustness of the DSMC is guaranteed, but also the conservative controller design, which occurs by using the maximum uncertainty bounds, is avoided. 

{In addition to implementation imprecisions issue, it has been shown in the literature that the SMC structure allows for handling the uncertainties in the plant model.}  The results in~\cite{Misawa_DSC} showed that the bounded uncertainties can be addressed without any adaptation by applying a saturation function into the SMC formulation that generates the boundary layer. The SMC with adaptation for handling modeling uncertainty/mismatch in the previous studies are limited to the continuous-time domain~\cite{Slotine}, linear systems~\cite{Chan_Automatica}, and first order sliding control~\cite{Pan_DSC,Amini_CEP,Amini_DSCC2016}. 

In this paper, a new adaptive second order DSMC formulation is developed for a general class of \textit{uncertain nonlinear systems}. Moreover, the asymptotic stability of the new controller is guaranteed via a Lyapunov stability argument and invoking the new Invariance Principle for nonlinear systems with discontinuity. {Compared to a conventional first order DSMC, the proposed adaptive second order DSMC with predicted implementation imprecisions provides (i) higher robustness against implementation imprecisions and modeling uncertainties, and (ii) faster tracking performance under unknown external disturbances. These benefits come at the cost of slightly more complex control logic (i.e., first order sliding mode versus second order sliding mode). However, the real-time test results show that the required computational power for the second order DSMC is almost the same as the first order DSMC. Thus, the real-time implementation of the proposed controller on a real ECU is feasible and it does not add any further computational demand, compared to the first order DSMC.} 

The proposed adaptive second order DSMC is generic and can be applied to broad engineering system applications. However, the most notable impact is anticipated for industry applications with intensive verification and validation (V\&V) practices, like automotive industries. This is because imprecisions which arise during controller implementation, and uncertainties within the model lead to iterative and costly V\&V. In this paper, application of the proposed control design is tested for control of a spark ignition (SI) combustion engine which exhibits highly nonlinear and coupled dynamics, and includes ADC and modeling uncertainties with required processing time of 2-10 \textit{ms} per iteration. {The contribution of this paper is threefold: \vspace{-0.1cm}
\begin{enumerate}
\item{Developing a new formulation of an \textit{adaptive second order DSMC} for a general class of MIMO nonlinear affine systems with inclusion of uncertainties within the plant model.}
\item Derivation of an innovative adaptation law to remove the uncertainty within the model in finite time via a discrete Lyapunov stability analysis, that also guarantees the asymptotic stability of the closed-loop system.
\item Presenting the first application of the adaptive second order DSMC for a combustion engine control problem. The proposed second order DSMC not only demonstrates robust behavior against data sampling and quantization imprecisions, but also removes the uncertainties in the engine model quickly and steers the dynamics to their nominal values in finite time.
\end{enumerate}
}
\vspace{-0.65cm}
\section{Second Order Discrete Sliding Mode Control} \label{sec:SISO2DSMC}

Let us consider the nonlinear system defined by:
\vspace{-0.45cm}
\begin{gather}\label{eq:C2SMC_1}
\dot{\mathbf{x}}(t)=f(t,\mathbf{x},u)
\end{gather}
where $\mathbf{x}=[x_1,~x_2,\cdots,x_r]^{\intercal}\in{\mathbb{R}^{r}}$, and $u \in{\mathbb{R}}$ are the state vector, and the scalar input variable, respectively. {Moreover, it is assumed that $f$ is smooth and sufficiently differentiable~\cite{salgado2004robust}.} The sliding mode order is the number of continuous successive derivatives of the differentiable sliding variable ${s}\in{\mathbb{R}}$, and it is a measure of the degree of smoothness of the sliding variable in the vicinity of the sliding manifold.

{The affine single-input single-output (SISO) form of the nonlinear system in Eq.~(\ref{eq:C2SMC_1}) with an unknown multiplicative term ($\alpha_i$) can be presented using the following state space equation:
\vspace{-0.55cm}
\begin{gather}\label{eq:D2SMC_1}
\dot{x}_i(t)=\alpha f_i(x_i(t))+g_i(x_i(t))u_i(t),
\end{gather}
where $g_i(x_i(t))$ is a non-zero input coefficient and $f_i(x_i(t))$ represents the dynamics of the plant and does not depend on the inputs. Note that the subscripts ($i$) are provided to represent a single scalar value. $\alpha_i$, which represents the errors/mismatches in the plant model, is unknown and constant. The continuous-time model in Eq.~(\ref{eq:D2SMC_1}) is discretized by utilizing the first order Euler approximation:
\vspace{-0.55cm}
\begin{gather}\label{eq:D2SMC_2}
x_i(k+1)=T\alpha_i f_i(x_i(k))+Tg_i(x_i(k))u_i(k)+x_i(k),
\end{gather}
in which $T$ is the sampling time. The second order discrete sliding variable is defined~\cite{mihoub2009real}:
\vspace{-0.65cm}
\begin{gather}\label{eq:D2SMC_3}
\xi_i(k)={s_i}(k+1)+\beta s_i(k),~\beta>0,
\end{gather}
where $s_i(k)=x_i(k)-x_{i,d}(k)$, and $\beta$ is a constant second order sliding mode gain. Moreover 
$x_{i,d}$ is the known desired trajectory of $x_{i}$. 

The equivalent control input of the second order DSMC should satisfy the following second order discrete sliding mode criteria~\cite{mihoub2009real}:
\vspace{-0.65cm}
\begin{gather}\label{eq:D2SMC_5}
\xi_i(k+1)=\xi_i(k)=0.
\end{gather}
{Applying Eq.~(\ref{eq:D2SMC_5}) to the nonlinear system in Eq.~(\ref{eq:D2SMC_2}) yields:
\vspace{-0.45cm}
\begin{gather}\label{eq:D2SMC_5_New1}
\xi_i(k)=0: x_i(k+1)-x_{i,d}(k+1)+\beta(x_i(k)-x_{i,d}(k))=0
\end{gather}
By substituting Eq.~(\ref{eq:D2SMC_2}) in Eq.~(\ref{eq:D2SMC_5_New1}) we have:
\vspace{-0.45cm}
\begin{gather}\label{eq:D2SMC_5_New2}
T\alpha_if_i(x_i(k))+Tg_i(x_i(k))u_i(k)+x_i(k)-x_{i,d}(k+1)\\
+\beta(x_i(k)-x_{i,d}(k))=0 \nonumber
\end{gather}
Based on which $u_i(k)$ can be calculated. However, since the value of $\alpha_i$ is unknown, the equivalent control input ($u_{i,eq}$) is expressed as function of the unknown $\alpha_i$ ($\hat{\alpha}_i$)  as follows:} 
\vspace{-0.45cm}
\begin{gather}\label{eq:D2SMC_6}
u_{i,eq}(k)=
\frac{-1}{g_iT}\Big(T\hat{\alpha_i}(k)f_i(x_i(k))+x_i(k)
-x_{i,d}(k+1)+\beta s_{i}(k)\Big), 
\end{gather}
By incorporating the control law ($u_{i,eq}$) into the second order sliding variable ($\xi_i$), we have:
\vspace{-0.45cm}
\begin{gather}\label{eq:D2SMC_7}
\xi_i(k)=Tf_i(x_i(k))(\alpha_i-\hat{\alpha}_i(k))=Tf_i(x_i(k))\tilde{\alpha}_i(k).
\end{gather}
$\tilde{\alpha}_i(k)$ is the error in estimating the unknown multiplicative term ($\tilde{\alpha}_i(k)=\alpha_i-\hat{\alpha}_i(k)$). Next, a Lyapunov stability analysis is conducted to (i) determine the stability of the closed-loop system, and (ii) derive the adaptation law to remove the uncertainty in the model. To this end, the following Lyapunov candidate function is proposed: 
\vspace{-0.5cm}
\begin{gather}\label{eq:D2SMC_8}
V_i(k)=\frac{1}{2}\Big({s}_i^2(k+1)+\beta{s}_i^2(k)\Big)+\frac{1}{2}\rho_{\alpha}\Big(\tilde{\alpha}_i^2(k+1)+\beta \tilde{\alpha}_i^2(k)\Big), 
\end{gather}
where $\rho_{\alpha}>0$ is a tunable adaptation gain chosen for the numerical sensitivity of the unknown parameter estimation. As can be seen from Eq.~(\ref{eq:D2SMC_8}), the proposed Lyapunov function is positive definite and quadratic with respect to the sliding variable ($s_i(k)$) and the unknown parameter estimation error ($\tilde{\alpha}_i(k)$). The desired condition is asymptotic and finite-time convergence of both $s_i$ and $\tilde{\alpha}_i$ to zero, {which is guaranteed by exploiting the following results.} \vspace{0.1cm}

{$\blacksquare$ \textbf{Proposition~1}: ~\textit{If} the following adaptation law is used to estimate the unknown multiplicative term ($\alpha_i$):\vspace{-0.45cm}
\begin{gather}\label{eq:D2SMC_14}
\tilde{\alpha}_i(k+1)=\tilde{\alpha}_i(k)-\frac{Ts_i(k)f_i(x_i(k))}{\rho_{\alpha}}. 
\end{gather}
\textit{then} the difference function of the proposed Lyapunov function in Eq.~(\ref{eq:D2SMC_8}) becomes: \vspace{-0.5cm}
\begin{gather}\label{eq:D2SMC_15}
\Delta V_i(k)=-(\beta+1)\big(s_i^2(k+1)+\beta s_i^2(k)\big)\\
+O\left(\Delta {s}_i^2(k),\Delta{s}_i^2(k+1),\Delta\tilde{\alpha}_i^2(k),\Delta\tilde{\alpha}_i^2(k+1)\right). \nonumber
\end{gather} 
where, $\Delta s_i(k)\equiv s_i(k+1)-s_i(k)$ and $\Delta \tilde{\alpha}_i(k) \equiv \tilde{\alpha}_i(k+1)-\tilde{\alpha}_i(k)$. {The proof of Proposition~1 is presented in the Appendix section of the supplementary material document of this paper presented in~\cite{Amini_JDSMC2018_sup}}. $\blacksquare$ } \vspace{0.1cm}

Next, the proposed equivalent control input (Eq.~(\ref{eq:D2SMC_6})) of the second order DSMC is evaluated on Gao's reaching law~\cite{gao1995discrete}. To this end, we begin with the first order continuous-time SMC, which can be realized by applying the following sliding reaching law~\cite{Slotine,Xiao_HOSMC}:\vspace{-0.45cm}
\begin{gather}\label{eq:CST_1}
\dot{s}_i(t)=-\lambda s_i(t)-\varepsilon. sgn(s_i(t))
,~~~~\lambda>0,~\varepsilon>0.
\end{gather}
Eq.~(\ref{eq:CST_1}) can be discretized as follows:
\vspace{-0.45cm}
\begin{gather}\label{eq:CST_2}
s_i(k+1)=(1-T\lambda)s_i(k)-T\varepsilon. sgn(s_i(k)), 
\end{gather}
where the $1>1-T\lambda>0$ condition should be met to guarantee that states of the system will move monotonically toward the switching plane ($s_i=0$), and cross it in finite time~\cite{gao1995discrete}. Next, the first order sliding variable ($s_i$) in Eq.~(\ref{eq:CST_2}) is replaced with the second order sliding variable ($\xi_i$):
\vspace{-0.45cm}
\begin{gather}\label{eq:CST_3}
\xi_i(k+1)=(1-T\lambda)\xi_i(k)-T\varepsilon. sgn(\xi_i(k)).
\end{gather}
{$\blacksquare$ \textbf{Proposition~2}: ~\textit{If} the second order sliding function ($\xi_i$) has the dynamic shown in Eq.~(\ref{eq:CST_3}),~\textit{then} by applying the calculated equivalent control input of the second order DSMC from Eq.~(\ref{eq:D2SMC_6}) along with the so-called switching control input ($u_{i,sw}$), it can be concluded that:
\vspace{-0.45cm}
\begin{gather}\label{eq:CST_13_3}
{s}_i(k+2)\approx\beta^2s_i(k),~~\beta>0.
\end{gather}
{The proof of Proposition~2 is presented in the supplementary material document~\cite{Amini_JDSMC2018_sup}}.~$\blacksquare$ } \vspace{0.1cm}

Proposition~2 (Eq.~(\ref{eq:CST_13_3})) states that, if $1>\beta>\beta^2>0$, the proposed equivalent control input of the second order DMSC along with the switching function in Eq.~({B.9}) of~\cite{Amini_JDSMC2018_sup} fulfills Gao's discrete sliding mode reaching law, which guarantees finite time convergence of the system's states to the sliding manifold~\cite{gao1995discrete}.

In the next step, by expanding the second order terms ($O(.)$) in Eq.~(\ref{eq:D2SMC_15}), and assuming a small enough sampling time ($T$), which means all terms that contain $T^2$ can be neglected, Eq.~(\ref{eq:D2SMC_15}) can be re-arranged as follows: \vspace{-0.45cm}
\begin{gather}\label{eq:D2SMC_16}
\Delta V_i(k)=-\frac{1}{2}\Big(\beta s_i^2(k+1)+\beta s_i^2(k)+2\beta^2 s_i^2(k)\\
+s_i^2(k+1)-s_i^2(k+2)\Big) \nonumber \\
-\beta s_i(k+1)s_i(k)-s_i(k+2)s_i(k+1), \nonumber
\end{gather}
in which, it is assumed that the uncertainty in the model is compensated by applying Eq.~(\ref{eq:D2SMC_14}). According to the second order sliding variable $\xi_i$ definition (Eq.~(\ref{eq:D2SMC_5})) and Eq.~(\ref{eq:CST_13_3}), $s_i(k+1)$ and $s_i(k+2)$ can be replaced by $-\beta s_i(k)$ and $\beta^2 s_i(k)$, respectively. By doing these replacements, Eq.~(\ref{eq:D2SMC_16}) can be simplified as: \vspace{-0.45cm}
\begin{gather}\label{eq:D2SMC_18}
\Delta V_i(k)=-\frac{1}{2}\beta\big(-{\beta}^3-{\beta}^2+\beta+1\big)s_i^2(k).
\end{gather}
$-{\beta}^3-{\beta}^2+\beta+1$ is positive if $1>\beta>0$. In other words, if $1>\beta>0$, then $\Delta V_i(k)\leq 0$, which guarantees the stability of the system ($\xi_i(k) \rightarrow 0$). Of interesting note is the $1>\beta>0$ constraint which is consistent with the earlier concluded condition on $\beta$ to make sure Gao's sliding reaching law is fulfilled by the equivalent control input of the second order DSMC. 

We showed that by using $u_{eq}$ from Eq.~(\ref{eq:D2SMC_6}), $u_{i,sw}$ from Eq.~(B.9) of~\cite{Amini_JDSMC2018_sup}, and the adaptation law from Eq.~(\ref{eq:D2SMC_14}), the negative semi-definite condition of the Lyapunov difference function is guaranteed. This means that the sliding variable (the tracking error, $s_i$) converges to zero, and the error in estimating the unknown parameter ($\tilde{\alpha}_i$) is at least bounded. Moreover, according to Eq.~(\ref{eq:CST_13_3}), since $1>\beta>0$, it is obvious that $1>\beta>\beta^2>0$. {Existence of a negative semi-definite Lyapunov difference function means that:
\vspace{-0.45cm}
\begin{gather}\label{eq:CST_14}
s_i(k)\geq s_i(k+1)\geq s_i(k+2),
\end{gather}}
{On the other side, based on the results of Proposition 2, we have ${s}_i(k+2)\approx\beta^2s_i(k)$. Thus, $s_i(k+2)=s_i(k)$, only if $s_i(k+2)=s_i(k)=0$; otherwise, $s_i(k+2)<s_i(k)$. We also know that based on the second order sliding control law $s_i(k+2)+\beta s_i(k+1)=0$ and $s_i(k+1)+\beta s_i(k)=0$. Therefore, if $s_i(k+2)=s_i(k)=0$, then $s_i(k+1)=0$. Overall, based on the negative semi-definite Lyapunov difference function, it can be shown that $s_i(k+2)<s_i(k)$, unless $s_i(k)=s_i(k+1)=s_i(k+2)=0$.} {Therefore, the second order DSMC guarantees asymptotic decrease of $s$ over a two-step horizon.} However, in order to guarantee the global asymptotic convergence and stability characteristics of the second order DMSC over the reaching and sliding modes, it is required to show that all the higher order Lyapunov difference functions vanish as $s_i\rightarrow 0$~\cite{Selina_PhD,Amini_CEP}. This will be proofed in the following section by invoking the new Invariance Principle for nonautonomous systems~\cite{barkana2015new}. \vspace{-0.45cm}

\subsection{Global Asymptotic Stability of the 2$^{nd}$ Order DSMC}

In order to guarantee the global asymptotic convergence and stability characteristics of the second order DSMC, the analysis begins with the second order Lyapunov difference function ($\Delta V_i(k+1)-\Delta V_i(k)$). Let us first define a new term, $\Gamma=\beta(-{\beta}^3-{\beta}^2+\beta+1)$, where $\Gamma>0$. With the definition of $\Gamma$, $\Delta V_i(k+1)-\Delta V_i(k)$ is calculated as follows:\vspace{-0.45cm}
\begin{gather}\label{eq:Model_un_newLya5}
\Delta V_i(k+1)-\Delta V_i(k)\approx \frac{1}{2}\Gamma{s}_i^2(k+1)-\frac{1}{2}\Gamma{s}_i^2(k) \\
\approx\frac{1}{2}\Gamma({s}_i^2(k+1)-{s}_i^2(k)). \nonumber
\end{gather}

The objective is to show that $\Delta V_i(k+1)-\Delta V_i(k)=0$ when $s_i(k)=0$. According to Eq.~(\ref{eq:Model_un_newLya5}), it is obvious that $\Delta V_i(k+1)-\Delta V_i(k)=0$ if ${s}_i^2(k+1)-{s}_i^2(k)=0$. ${s}_i^2(k+1)-{s}_i^2(k)=0$ if either of the following conditions is met:\vspace{-0.45cm}
\begin{subequations}\label{eq:Model_un_newLya6}
\begin{gather}
s_i(k+1)=s_i(k)=0~~~(I) \\
s_i(k+1)=s_i(k)\neq 0~~~(II).
\end{gather}
\end{subequations}
Eq.~(\ref{eq:D2SMC_18}) says that if $s_i(k)=0$, then $\Delta V_i(k)=0$. If one assumes that the unknown uncertainty term is removed from the model using the proposed adaptation algorithm, the second order sliding mode condition denotes that $\xi_i(k+1)=\xi_i(k)=0$, which means $s_i(k+1)=-\beta s_i(k)$. Thus, if $s_i(k)=0$, $s_i(k+1)$ also becomes zero, and consequently $\Delta V_i(k+1)-\Delta V_i(k)=0$. This means that condition (I) in Eq.~(\ref{eq:Model_un_newLya6}) is realizable.

According to Eq.~(\ref{eq:CST_14}), $s_i(k+1)< s_i(k)$. This means $s_i(k+1)\neq s_i(k)$ unless $s_i(k+1)=s_i(k)=0$. Thus, the condition ($II$) in Eq.~(\ref{eq:Model_un_newLya6}) cannot be true and only condition ($I$) is feasible. According to condition ($I$), if $s(k)=0$, not only $\Delta V_i(k)=0$, but also $\Delta V_i(k+1)-\Delta V_i(k)=0$. 

Next, the third Lyapunov difference function should be calculated:\vspace{-0.45cm}
\begin{gather}\label{eq:Model_un_newLya7}
[\Delta V_i(k+2)-\Delta V_i(k+1)]-[\Delta V_i(k+1)-\Delta V_i(k)] \\
\approx\frac{1}{2}\Gamma({s}_i^2(k+2)-2{s}_i^2(k+1)+{s}_i^2(k)). \nonumber 
\end{gather}
Again, since $s_i(k+2)< s_i(k+1)<s_i(k)$, $s_i(k+2)\neq s_i(k+1)\neq s_i(k)$ unless $s_i(k+2)=s_i(k+1)=s_i(k)=0$. Thus, if $s_i(k)=0$, not only $s_i(k+1)=s_i(k+2)=0$, but also $\Delta V_i(k)=0$, $\Delta V_i(k+1)-\Delta V_i(k)=0$, and $(\Delta V_i(k+2)-\Delta V_i(k+1))-(\Delta V_i(k+1)-\Delta V_i(k))=0$. In a same manner as the first, second, and third order difference functions, it can be shown that higher order Lyapunov difference functions ($>$3) become zero only, and only if $s_i(k)=0$~\cite{Selina_PhD,Amini_CEP}. This is a key conclusion that allows for proof of the global asymptotic stability by invoking the new Invariance Principle for nonautonomous systems~\cite{barkana2014defending}.

Continuity is one of the required conditions for the LaSalle's Invariance Principle~\cite{khalil1996nonlinear} to conclude the asymptotic stability with respect to a negative semi-definite derivative of a positive definite Lyapunov function. LaSalle's Invariance Principle has been extended in~\cite{barkana2014defending} to nonlinear systems with discontinuity. The extension of LaSalle's Invariance Principle to discrete systems, which is called the new Invariance Principle theorem, removes the continuity requirement, and allows us to conclude the asymptotic stability with respect to a negative semi definite difference function of a positive definite Lyapunov discrete equation. It was shown that when $s_i(k)=0$, the Lyapunov difference function, all the future values of the sliding variable, and higher order Lyapunov differences become zero. Therefore, all the trajectories of the system approach the set defined by $\Delta V_i(k)\equiv0$. Since $s(k)=0$ and $\tilde{\alpha}_i(k)=0$ are the only trajectories which satisfy the nonlinear uncertain system equations, this trajectory ($s(k)=0,~\tilde{\alpha}_i(k)=0$) is a limit point, and also an equilibrium point, of the closed-loop system~\cite{Selina_PhD,Amini_CEP}.

Since the Lyapunov difference function cannot be negative for any unlimited period of time, according to the new Invariance Principle theorem~\cite{barkana2014defending}, $\Delta{V}_i$ must be identically zero at any limit point~\cite{barkana2015new}. Let define two new domains:\vspace{-0.45cm}
\begin{gather}
\label{eq:Model_un_newLya8}
\Omega_0=\{x_i|V_i(x_i)\leq V(x_{i,0})\} \\
\Omega_i=\{x_i|\Delta V_i(x)\equiv0)\}, \nonumber
\end{gather}
where $x_{i,0}$=$x_i(0)$ is the initial condition. The negative semi definite condition of the Lyapunov difference function means that all system states are bounded and contained within the domain $\Omega_0$~\cite{barkana2014defending}. For a small enough sampling time, if the following condition holds~\cite{barkana2014defending}:

\vspace{0.1cm}
\hspace{-0.2cm}$|x_i(k)+Tf_i+Tg_iu_i(k)|$ \textit{is bounded for any bounded $x_i$}
\vspace{0.1cm}

then, $s_i(k)$=$\tilde{\alpha}_i(k)$=$0$, which is the only limit point (equilibrium point) of the system, belongs to $\Omega_f$=$\Omega_0 \bigcap \Omega_i$. Therefore, $s_i(k)$ and $\tilde{\alpha}_i(k)\rightarrow 0$ as $k\rightarrow\infty$. Thus, according to the new Invariance Principle theorem~\cite{barkana2014defending}, the asymptotic convergence of the two variables to zero with a positive definite Lyapunov function and a negative semi-definite Lyapunov difference equation is concluded.

Overall, the control input of the second order DSMC is:
\vspace{-0.65cm}
\begin{gather}\label{eq:CST_21}
u_i(k)=u_{i,eq}(k)+u_{i,sw}(k),
\end{gather}
where $u_{i,eq}$ and $u_{i,sw}$ are calculated according to Eq.~(\ref{eq:D2SMC_6}), and Eq.~(B.9) of~\cite{Amini_JDSMC2018_sup}, respectively. The main reason to add the switching function to the calculated equivalent control input of the SMC is reducing the chattering~\cite{nollet2008observer}. As it has been shown in~\cite{Amini_CEP,Kyle_DSCC,Amini_ACC2016}, the addition of the switching function with a fix or variable gain enhances the robustness of the controller against uncertainties, which could overpower the chattering issues. The gain of the switching function ($sgn(\xi_i(k-1))$), which depends on the DSMC tuning parameters ($\lambda,~\beta$) and $g$ (Eq.~(\ref{eq:D2SMC_1})), is the key to ensure the controller robustness against external uncertainties~\cite{Amini_CEP}. By looking into Eq.~(B.9) of~\cite{Amini_JDSMC2018_sup}, it can be seen that the gain of the switching function has the same unit of the control input. This gain represents the boundary of the external uncertainties on the control signals, i.e. ADC imprecisions~\cite{Kyle_DSCC,Amini_DSC}. 

Tuning the switching function gain based on constant $\lambda$ and $\beta$ is hard to achieve. Instead, an online sampling and quantization uncertainty prediction mechanism was proposed in~\cite{Amini_ACC2016} which allows for estimating and propagating the ADC uncertainty bounds on the control signal, and avoids conservative controller design. The switching function with the predicted ADC uncertainty bounds has the following structure: 
\vspace{-0.65cm}
\begin{gather}\label{eq:CST_22}
u_{i,sw}(k)=-|\mu_{u_i}|sgn(\xi_i(k-1)), 
\end{gather}
where $\mu_{u_i}$ is the propagated sampling and quantization imprecisions. $\mu_{u_i}$ is estimated according to the mechanism presented with details in~\cite{Amini_CEP,Amini_ACC2016}. Fig.~\ref{fig:AdaptiveDSMC_Schematic} illustrates the overall schematic of the proposed adaptive second order DSMC along with the ADC uncertainty and propagation mechanism. In order to avoid possible high frequency chattering which occurs in discrete systems during implementation of the signum function, the signum function is replaced with saturation ($sat$) function which provides smoother behavior:
\vspace{-0.65cm}
\begin{gather}\label{eq:CST_23}
u_{i,sw}(k)=-|\mu_{u_i}|sat\Big(s_i(k)+\beta s_i(k-1)\Big). 
\end{gather}
\vspace{-0.6cm} 
\begin{figure}[h!]
\begin{center}
\includegraphics[width=0.9\columnwidth]{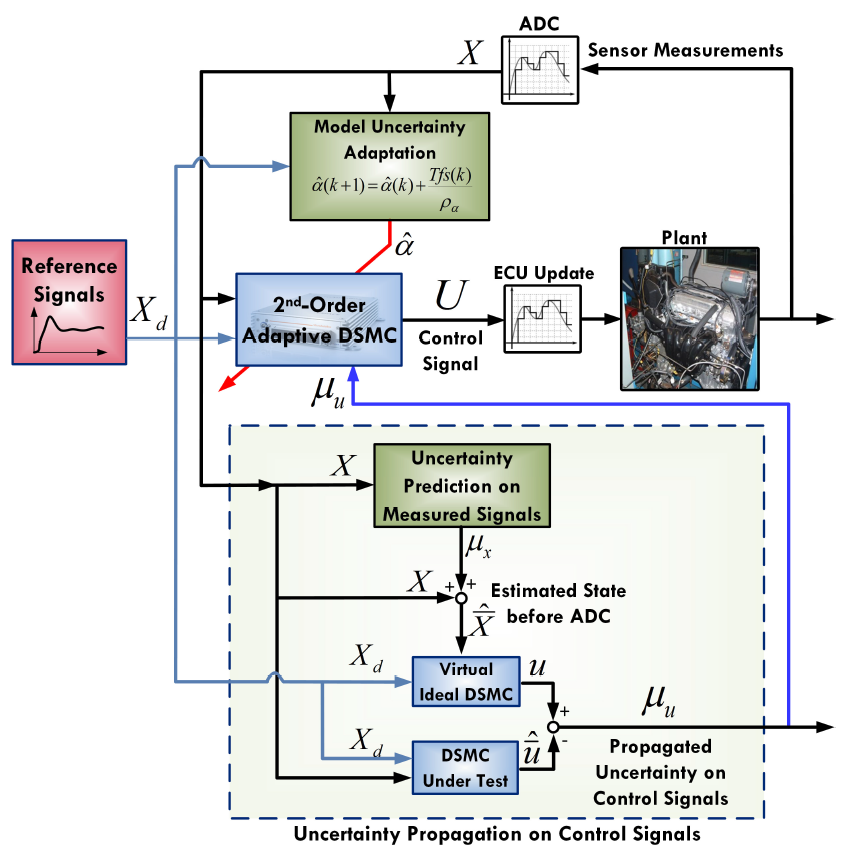} \vspace{-0.45cm}
\caption{\label{fig:AdaptiveDSMC_Schematic} Schematic of the second order adaptive DSMC with ADC uncertainty prediction and propagation mechanism.} \vspace{-0.8cm}
\end{center}
\end{figure}

\vspace{-0.65cm}
\section{MIMO Adaptive Second Order DSMC} 
\label{sec:MIMOSecondDSMC}
The discretized nonlinear system can be expressed in MIMO structure as follows:
\vspace{-0.45cm}
\begin{gather}\label{eq:MIMO_1}
\mathbf{x}(k+1)=T\mathbf{f}(\mathbf{x}(k))+T\mathbf{g}(\mathbf{x}(k))\mathbf{u}(k)+\mathbf{x}(k),
\end{gather}
where $x_i\subset{\mathbf{x}\in{\mathbb{R}^{r}}},~(i=1:r)$, and $u_j\subset{\mathbf{u}\in{\mathbb{R}^{h}}},~(j=1:h)$ are the state and control input vectors, respectively. The remaining dynamics are represented as $\mathbf{f}(\mathbf{x}(k))$ and $\mathbf{g}(\mathbf{x}(k))$. A generic first order DSMC for the MIMO system in Eq.~(\ref{eq:MIMO_1}) has been proposed in our previous work~\cite{Pan_DSC}. Here, we use the results from \cite{Pan_DSC} to derive the first order sliding surface vector ($\mathbf{s}(k)=[s_1(k),...,s_h(k)]^\intercal$) based on a system output vector which is defined to be $y_j\subset{\mathbf{y}\in{\mathbb{R}^{h}}},~(j=1:h)$:
\vspace{-0.45cm}
\begin{gather}\label{eq:MIMO_2}
y_j(k)=m_j(\mathbf{x}(k)).
\end{gather}

For a system with an output ${y}_j$ that has a relative degree of $\kappa_j$, the scalar sliding surface $s_j(t)$ is defined as:
\vspace{-0.45cm}
\begin{gather}\label{eq:MIMO_2_1}
s_j(t)=\Big(\frac{d}{dt}+\lambda_j\Big)^{\kappa_j-1}(y_j(t)-y_{j_{d}}(t)),
\end{gather}
where $y_{j_d}\subset{\mathbf{y}_d\in{\mathbb{R}^{h}}},~(j=1:h)$ is the desired output values. It was shown in~\cite{Pan_DSC} that for a MIMO system with a relative degree of $\kappa_j$ ($\sum_{1}^{h}\kappa_j<r$), the discrete dynamics of the first order sliding surface ($\mathbf{s}$) in the presence of the model uncertainties is: 
\vspace{-0.45cm}
\begin{gather}\label{eq:MIMO_4}
\mathbf{s}(k+1)=\mathbf{s}(k)+T\hat{\boldsymbol{\Lambda}}+T\mathbf{F}\mathbf{a}+T\boldsymbol{\Upsilon}\mathbf{u}(k),
\end{gather}
where $\mathbf{a}\in{\mathbb{R}^{p}}$ is a vector of unknown constants in the model and $\mathbf{F}\in{\mathbb{R}^{h\times p}}$ is known as a data matrix~\cite{Pan_DSC}. $\hat{\boldsymbol{\Lambda}}$ is defined as $\hat{\boldsymbol{\Lambda}}={\boldsymbol{\Lambda}}-\mathbf{F}\mathbf{a}$, where ${\boldsymbol{\Lambda}}=[l_1,l_2,...,l_h]^\intercal$, and:
\vspace{-0.45cm}
\begin{gather}\label{eq:MIMO_4_3}
l_j=L^{\kappa_j}_f(h_j)+\sum_{q=\kappa_j-1}^{1}c_{j(\kappa_j-q+1)}.[L^{q}_f(h_j)-\frac{d^{\kappa_j}}{dt^{\kappa_j}}y_{j_d}^{(\kappa_j)}(t)],
\end{gather}\vspace{-0.15cm}
and, $L^{q}_f(h_j)$ is defined as~\cite{khalil1996nonlinear}:
\vspace{-0.45cm}
{
\begin{gather}\label{eq:MIMO_4_4}
L^{q}_f(m_j)=\frac{d^q{y_j}(t)}{dt^q}, 
\end{gather}}
with $c_{j(\kappa_j-q+1)}$ chosen such that all poles are at $-\lambda_j$~\cite{Selina_PhD}. In the absence of the model uncertainties, it is obvious that $\hat{\boldsymbol{\Lambda}}={\boldsymbol{\Lambda}}$. In this paper, we assume a relative degree of one ($\kappa_j=1$) for the output. Thus, according to the sliding surface definition and the relative degree of the outputs, the non-singular $\boldsymbol{\Upsilon}\in{\mathbb{R}^{h\times h}}$ matrix in Eq.~(\ref{eq:MIMO_4}) becomes:
\vspace{-0.45cm}
\begin{gather}\label{eq:MIMO_5}
\boldsymbol{\Upsilon}=\begin{bmatrix}
g_{11} & ... & g_{1h}\\ 
\vdots  & \ddots  & \\ 
g_{h1} & ... & g_{hh}
\end{bmatrix}.
\end{gather}

It is assumed that here we deal with a MIMO system that has an output with a relative degree of one. Therefore, according to Eq.~(\ref{eq:MIMO_2_1}), $s_j(k)$=$y_j(k)-{y_d}_j(k)$. For the sake of simplicity it can be assumed that $r=h$ (which means $i$=$j$=$p$=1:$r$), $m_i(\mathbf{x})$=$x_i(k)$, and ${y_d}_i(k)$=${x_d}_i(k)$. The latter assumption means that here, state space variables of the state vector are treated as the output variables, and for each state variable ($x_i$), a sliding surface ($s_i$) is defined and it is assumed that a unique control input ($u_i$), either physical or synthetic, exists for every single sliding surface within the sliding surface vector ($\mathbf{s}$). 

Now, the second order sliding surface vector ($\boldsymbol{\xi}=[\xi_1,...,\xi_r]^\intercal$) is constructed:
\vspace{-0.45cm}
\begin{gather}\label{eq:MIMO_6}
\boldsymbol{\xi}(k)=\mathbf{s}(k+1)+\boldsymbol{\beta}\mathbf{s}(k),
\end{gather}
where $\mathbf{s}(k+1)$ is calculated with respect to Eq.~(\ref{eq:MIMO_4}) and $\boldsymbol{\beta}\in{\mathbb{R}^{r\times r}}$ is the \textit{positive definite} matrix of the second order sliding mode gains. The equivalent control input vector which satisfies the second order sliding mode condition ($\boldsymbol{\xi}(k+1)=\boldsymbol{\xi}(k)=0$) can be found by solving the following equation for $\mathbf{u}_{eq}$:
\vspace{-0.45cm}
\begin{gather}\label{eq:MIMO_7}
T\boldsymbol{\Upsilon}\mathbf{u}_{eq}(k)=-(\mathbf{I+\boldsymbol{\beta}})\mathbf{s}(k)-T\hat{{\boldsymbol{\Lambda}}}-T\mathbf{F}\hat{\mathbf{a}}(k),
\end{gather}
where $\hat{\mathbf{a}}(k)$ is the estimation of the unknown constants ${\mathbf{a}}$. In a similar manner to the SISO controller, upon substitution of Eq.~(\ref{eq:MIMO_7}) into Eq.~(\ref{eq:MIMO_6}), we have:
\vspace{-0.45cm}
\begin{gather}\label{eq:MIMO_8}
\boldsymbol{\xi}(k)=T\mathbf{F}\tilde{\mathbf{a}}(k),
\end{gather}
where $\tilde{\mathbf{a}}(k)={\mathbf{a}}-\hat{\mathbf{a}}(k)$. 

Similar to the SISO controller, a discrete Lyapunov analysis is carried out to first converge $\tilde{\mathbf{a}}$ to zero and remove the uncertainty in the model, and second, guarantee the asymptotic stability of the closed loop system. To this end, the argument begins with the following Lyapunov candidate function~\cite{Pan_DSC}:
\vspace{-0.45cm}
\begin{gather}\label{eq:MIMO_9_a2}
{V}(k)=\frac{1}{2}\mathbf{s}(k)^\intercal\mathbf{s}(k)+\frac{1}{2}\tilde{\mathbf{a}}(k)^\intercal\boldsymbol{\Gamma}\tilde{\mathbf{a}}(k),
\end{gather}
where $\boldsymbol{\Gamma}\in{\mathbb{R}^{r\times r}}$ is the tunable adaptation \textit{symmetric positive} matrix. \vspace{0.15cm}

{$\blacksquare$ \textbf{Proposition~3}: ~\textit{If} the following adaptation law is used to estimate the vector of the unknown multiplicative terms ($\tilde{\mathbf{a}}$):\vspace{-0.45cm}
\begin{gather} \label{eq:MIMO_17_a2}
\hat{\mathbf{a}}(k+1)=\hat{\mathbf{a}}(k)+T(\boldsymbol{\Gamma}\boldsymbol{\Gamma})^{-1}\mathbf{F}^\intercal\mathbf{s}(k),
\end{gather}
\textit{then} the Lyapunov difference function of Eq.~(\ref{eq:MIMO_9_a2}) becomes:
\vspace{-0.45cm}
\begin{gather} \label{eq:MIMO_18_2_a2}
\Delta{V}(k)=-\frac{1}{2}\mathbf{s}(k)^\intercal(\mathbf{I}-\boldsymbol{\beta}\boldsymbol{\beta})\mathbf{s}(k).
\end{gather}
where, $\Delta{V}$ is negative semi-definite if eigenvalues of $\boldsymbol{\beta}$ lie within the unit circle~\cite{Selina_PhD}. The proof of Proposition~3 is presented in the Appendix section of the supplementary material document~\cite{Amini_JDSMC2018_sup}. $\blacksquare$ } \vspace{0.45cm}

One can easily find the analogous structure between MIMO (Eq.~(\ref{eq:MIMO_17_a2})) and SISO (Eq.~(\ref{eq:D2SMC_14})) adaptation laws.

According to Proposition~3, the Lyapunov function in Eq.~(\ref{eq:MIMO_9_a2}) guarantees the finite-time zero convergence of the first order sliding vector ($\mathbf{s}$), and intuitively gives the adaptation law to remove the uncertainty in the model. However, still it is required to prove the asymptotic convergence of the second order sliding vector ($\boldsymbol{\xi}$) to zero. To this end, first it is assumed that the uncertainty in the model is removed permanently by incorporating the adaptation rule from Eq.~(\ref{eq:MIMO_17_a2}). Validity of this assumption will be testified in Sec.~\ref{sec:Results}. Next, the system switches to another Lyapunov function ($V^*$, proposed in Eq.~(\ref{eq:MIMO_new_9})), which covers the time interval after the completion of the adaptation period. 

We define the new Lyapunov function ($V^*$) as follows:
\vspace{-0.45cm}
\begin{gather}\label{eq:MIMO_new_9}
{V^*}(k)=\frac{1}{2}\Big(\mathbf{s}(k+1)^\intercal\mathbf{s}(k+1)+\mathbf{s}(k)^\intercal\boldsymbol{\beta}\mathbf{s}(k)\Big).
\end{gather}
By using the Taylor series expansion, implementing the values for the first and second order partial derivatives of $V^*$ with respect to $\mathbf{s}(k)$ and $\mathbf{s}(k+1)$, and knowing that all the higher order partial derivatives and also second order cross derivative are zero, the Lyapunov difference function becomes:
\vspace{-0.65cm}
\begin{gather}\label{eq:MIMO_new_10}
\Delta V^*(k)=\mathbf{s}(k)^\intercal\boldsymbol{\beta}(\mathbf{s}(k+1)-\mathbf{s}(k))\\
~~~~~~~~~~~~~~~~+\mathbf{s}(k+1)^\intercal(\mathbf{s}(k+2)-\mathbf{s}(k+1))\nonumber \\
~+\frac{1}{2}\Big(\Delta \mathbf{s}(k)^\intercal\boldsymbol{\beta}\Delta \mathbf{s}(k)+\Delta \mathbf{s}(k+1)^\intercal\Delta \mathbf{s}(k+1)\Big).\nonumber
\end{gather}
According to the earlier assumption, upon removal of the uncertainties in the model, Eq.~(\ref{eq:MIMO_8}) becomes $\boldsymbol{\xi}(k)=\mathbf{0}$. Thus, Eq.~(\ref{eq:MIMO_new_10}) can be simplified as follows:
\vspace{-0.45cm}
\begin{gather}\label{eq:MIMO_new_11}
\Delta V^*(k)=-\mathbf{s}(k)^\intercal\boldsymbol{\beta}(\boldsymbol{\beta}+\mathbf{I})\mathbf{s}(k)\\
~~~~~~~~~~~~~~~~-\mathbf{s}(k+1)^\intercal(\boldsymbol{\beta}+\mathbf{I})\mathbf{s}(k+1)\nonumber \\
~~~~~~~~+\frac{1}{2}\Big(\Delta \mathbf{s}(k)^\intercal\boldsymbol{\beta}\Delta \mathbf{s}(k)+\Delta \mathbf{s}(k+1)^\intercal\Delta \mathbf{s}(k+1)\Big).\nonumber
\end{gather}

In a similar manner to the SISO system, it can be shown that upon applying the equivalent control input of the MIMO second order DSMC (Eq.~(\ref{eq:MIMO_7})) along with the vector of the switching functions ($|\mu_{u_i}|sat(\xi_i)$) on the MIMO form of Gao's reaching law~\cite{Monsees_PhD}, and assuming the removal of the model uncertainties by using Eq.~(\ref{eq:MIMO_17_a2}), we obtain:
\vspace{-0.45cm}
\begin{gather} \label{eq:MIMO_19}
\mathbf{s}(k+2)\approx\boldsymbol{\beta}^2\mathbf{s}(k),
\end{gather}
where the input of the second order MIMO DSMC is~\cite{Monsees_PhD}:
\vspace{-0.45cm}
\begin{gather} \label{eq:MIMO_19_2}
\mathbf{u}=\mathbf{u}_{eq}(k)-
\begin{bmatrix}
|\mu_{u_1}|sat(\xi_1(k-1)) \\ 
\vdots  \\ 
|\mu_{u_r}|sat(\xi_i(k-1))
\end{bmatrix},
\end{gather}
in which $\mu_{u_i},~i=0:r$ is the predicted uncertainties on the corresponding control signal. Next, by substituting Eq.~(\ref{eq:MIMO_19}) into Eq.~(\ref{eq:MIMO_new_11}), we have:
\vspace{-0.45cm}
\begin{gather} \label{eq:MIMO_20}
\Delta{V^*}(k)=-\mathbf{s}^\intercal\boldsymbol{\beta}\Big(-\boldsymbol{\beta}^3-\boldsymbol{\beta}^2+\boldsymbol{\beta}-\mathbf{I}\Big)\mathbf{s}(k),
\end{gather}
which also has a similar structure to Eq.~(\ref{eq:D2SMC_18}). $\Delta{V^*}$ in Eq.~(\ref{eq:MIMO_20}), is negative semi-definite in $\mathbf{s}(k)$ when $\boldsymbol{\beta}$ is chosen such that its eigenvalues lie within the
unit circle. More importantly, due to the analogous structure between MIMO and SISO adaptive second order DSMCs, the asymptotic stability of the MIMO controller can be proved by invoking the new Invariance Principle~\cite{barkana2014defending,barkana2015new} and noting that the higher order Lyapunov difference functions are zero, if and only if $\mathbf{s}=\mathbf{0}$~\cite{Pan_DSC}.\vspace{-0.5cm} 

\section{Case Study: Automotive Engine Control}\label{sec:Engine_Control}
In this section, the proposed adaptive second order DSMC is designed for an experimentally validated physics-based SI engine model~\cite{Shaw} during transient cold start period. The engine model~\cite{Shaw} is parameterized for a 2.4-liter, 4-cylinder, DOHC 16-valve Toyota 2AZ-FE engine. The engine rated power is 117kW $@$ 5600~RPM, and it has a rated torque of 220 Nm $@$ 4000~RPM. The experimental validation of different components of the engine model is available in~\cite{Sanketi}. The nonlinear model has four states including the exhaust gas temperature~($T_{exh}$), fuel mass flow rate into the cylinders (${\dot{m}_f}$), the engine speed~(${\omega_e}$), and the mass of air inside the intake manifold ($m_{a}$). The engine performance is controlled by three inputs: (i) $\dot{m}_{ai}$ (air mass flow rate into the intake manifold) controls the engine speed via the air throttle body, (ii) $\dot{m}_{fc}$ (amount of injected fuel into the cylinder) regulates the Air-Fuel ratio ($AFR$) via the fuel injector, and (iii) $\Delta$ (spark timing) controls the exhaust gas temperature via the spark plug. Details of the functions and constants in the engine model are found in~\cite{Shaw}. The schematic of the SI engine control system is shown in Fig.~\ref{fig:Engine_Schematic}.
\begin{figure}[h!]
\begin{center}
\includegraphics[width= \columnwidth]{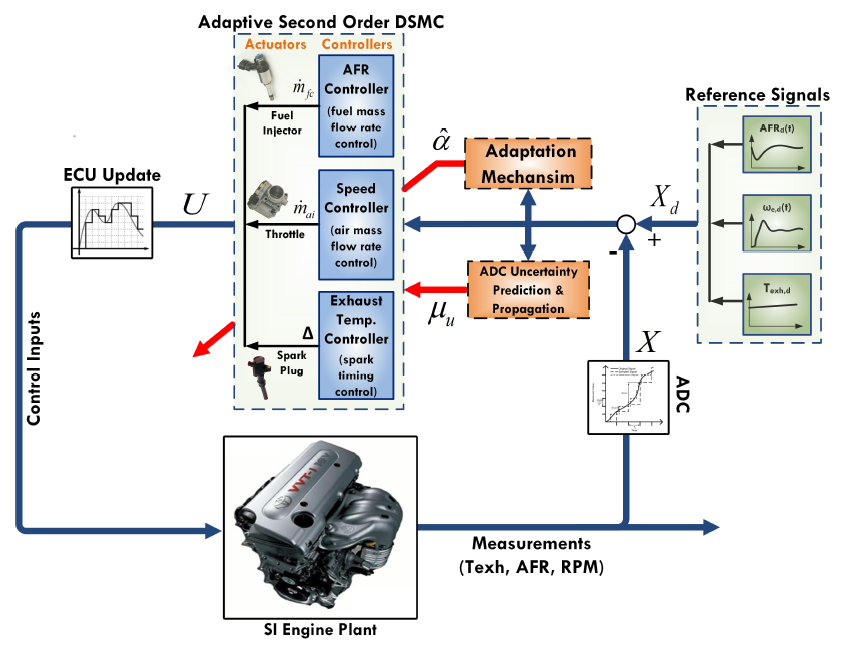} \vspace{-0.75cm}
\caption{\label{fig:Engine_Schematic} Block diagram of the engine cold start control system.} \vspace{-0.8cm}
\end{center}
\end{figure}

Four states of the model and corresponding nonlinear dynamics are as follows~\cite{Sanketi}:
\vspace{-0.45cm}
\begin{equation}\label{eq:Engine_IEEE_1}
\begin{split}
\begin{bmatrix}
T_{exh}(k+1)\\ 
\dot{m}_f(k+1)\\ 
m_a(k+1)\\ 
\omega_e(k+1)
\end{bmatrix}=
\begin{bmatrix}
T_{exh}(k)\\ 
\dot{m}_f(k)\\ 
m_a(k)\\ 
\omega_e(k)
\end{bmatrix}~~~~~~~~~~~~~~~~~~~~~~~~~~~~~~~~~\\
+T\left(\begin{bmatrix}
f_{T_{exh}}(k) & 0 & 0 & 0\\ 
0 & f_{\dot{m}_f}(k) & 0 & 0\\ 
0 & 0 & f_{m_a}(k) & 0\\ 
0 & 0 & 0 & f_{\omega_e}(k) 
\end{bmatrix}
\begin{bmatrix}
\alpha_{T_{exh}}\\ 
\alpha_{\dot{m}_f}\\ 
\alpha_{m_a}\\ 
\alpha_{\omega_e}
\end{bmatrix}\right.\\
\left.+\begin{bmatrix}
g_{T_{exh}}(k) & 0 & 0 & 0\\ 
0 & g_{\dot{m}_f}(k) & 0 & 0\\ 
0 & 0 & g_{m_a}(k) & 0\\ 
0 & 0 & 0 & g_{\omega_e}(k) 
\end{bmatrix}
\begin{bmatrix}
\Delta(k)\\ 
\dot{m}_{f,c}(k)\\ 
\dot{m}_{ai}(k)\\ 
m_{a,d}(k)
\end{bmatrix}
\right),	
\end{split}
\end{equation} \vspace{-0.25cm}
where $\mathbf{F}=diag[f_{{T_{exh}}},f_{\dot{m}_f},f_{\omega_e},f_{m_a}]$, and:
\begin{subequations}
\label{eq:Engine_IEEE_2}
\begin{align}
f_{{T_{exh}}}=\frac{1}{\tau_e}[600AFI-T_{exh}] 
\end{align}\vspace{-1.25cm}
\begin{align}
f_{\dot{m}_f}=-\frac{1}{\tau_f}\dot{m}_f(k)
\end{align}\vspace{-1.25cm}
\begin{align}
f_{\omega_e}=-\frac{1}{J}(0.4~\omega_e(k)+100)
\end{align}\vspace{-1.25cm}
\begin{align}
f_{m_a}=-\dot{m}_{ao}.
\end{align}
\end{subequations}
Also $\boldsymbol{\Upsilon}=diag[g_{{T_{exh}}},g_{\dot{m}_f},g_{\omega_e},g_{m_a}]$, and:\vspace{-0.45cm}
\begin{gather}
\label{eq:Engine_IEEE_3}
g_{{T_{exh}}}=\frac{7.5}{\tau_e},~~g_{\dot{m}_f}=\frac{1}{\tau_f},~~
g_{\omega_e}=\frac{30000}{J},~~g_{m_a}=1.
\end{gather}

The tracking control problem is defined to steer $T_{exh}$, ${\omega_e}$, $AFR$ to their desired values. To this end, with respect to each desired trajectory ($T_{exh,d}$, ${\omega_{e,d}}$, $AFR_d$), a sliding variable (tracking error) is defined. Additionally, as can be seen from Eq.~(\ref{eq:Engine_IEEE_1}), since there is no explicit control input in the engine speed equations, $m_{a,d}$ is defined as a synthetic control input for the rotational dynamics. The calculated $m_{a,d}$ will be used as the desired trajectory for the air mass flow controller over the next engine cycle. Thus, the sliding vector is:
\vspace{-0.45cm}
\begin{equation}\label{eq:Engine_IEEE_4}
\begin{split}
\mathbf{s}(k)=\begin{bmatrix}
s_1(k)\\ 
s_2(k)\\ 
s_3(k)\\ 
s_4(k)
\end{bmatrix}=
\begin{bmatrix}
T_{exh}(k)-T_{exh,d}(k)\\ 
AFR(k)-AFR_{d}(k)\\ 
m_a(k)-m_{a,d}(k)\\ 
\omega_e(k)-\omega_{e,d}(k)
\end{bmatrix}.
\end{split}
\end{equation}

The equivalent control input vector ($\mathbf{u}_{eq}$) of the $\textit{baseline}$ second order DSMC for the engine case study can be obtained by substituting $\mathbf{F}$, $\boldsymbol{\Upsilon}$, and $\mathbf{s}$ into Eq.~(\ref{eq:MIMO_7}). For the SISO DSMC, $\boldsymbol{\beta}$ is chosen to be diagonal, while for the MIMO controller, the dynamic coupling is included via the off-diagonal element of $\boldsymbol{\beta}$. For the second order DSMC \textit{with ADC uncertainties}, according to Eq.~(\ref{eq:MIMO_19_2}), the predicted ADC uncertainties are incorporated into the DSMC structure through the switching control input ($\mathbf{u}_{sw}$) gains. For the engine case study, the gains of the switching control input ($\mu_{\Delta}$, $\mu_{\dot{m}_{fc}}$, $\mu_{\dot{m}_{ai}}$, $\mu_{m_{a,d}}$) are estimated online using the mechanism previously shown in Fig.~\ref{fig:AdaptiveDSMC_Schematic}. The accuracy of the sampling and quantization imprecisions prediction and propagation mechanism has been shown in our previous works~\cite{Amini_CEP,Amini_DSCC2016,Amini_ACC2016}. {The engine controller is tuned manually, by using the try and error method, to achieve the best performance. The performed simulation and real-time results, which will be presented later in this section, show that once the controller is tuned properly, there is almost no need to re-tune the controller, unless the sampling and quantization level changes, which in practice those are known and fixed. }\vspace{-0.35cm}

\subsection{Handling Implementation Imprecisions}
In order to demonstrate the robustness characteristics of the second order DSMC compared to the first order controller in handling ADC uncertainties, first we assume that the engine model is ideal and there is no uncertainty in the modeled dynamics ($\alpha_{T_{exh}}$=$\alpha_{\dot{m}_f}$=$\alpha_{\omega_e}$=$\alpha_{{m}_a}$=1). Fig.~\ref{fig:Eng_Sampling80_20ms_SISO_1st_2nd} shows the desired trajectories ($AFR$, $T_{exh}$, and engine speed) tracking results, using the first and second order DSMCs for sampling times of $20~ms$ and $80~ms$, and quantization level of 16-$bit$ and 10-$bit$, respectively. The mean tracking errors for both controllers are listed in Table~\ref{table:tracking_Results}. It can be observed from Fig.~\ref{fig:Eng_Sampling80_20ms_SISO_1st_2nd} and Table~\ref{table:tracking_Results} that, when the signals at the controller I/O are sampled every $20~ms$, both first and second order baseline DSMCs illustrate smooth and acceptable tracking performances, while the second order controller is more accurate by up to 67\% in terms of the tracking errors. 
\begin{figure}[t!]
\begin{center}
\includegraphics[angle=0,width= \columnwidth]{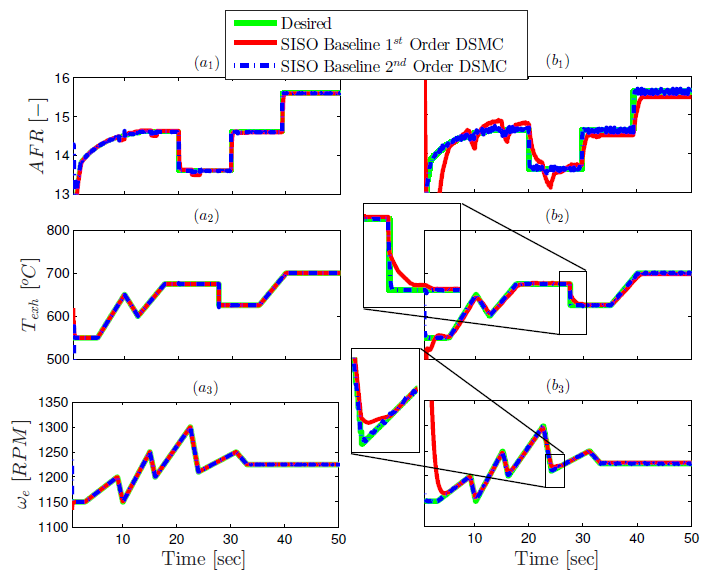} \vspace{-0.8cm}
\caption{\label{fig:Eng_Sampling80_20ms_SISO_1st_2nd}Engine tracking results by the first and second order SISO DSMCs for ($a$) $T$=20$~ms$, quantization level=16-$bit$, and ($b$) $T$=80$~ms$, quantization level=10-$bit$.} \vspace{-0.8cm}
\end{center}
\end{figure}
\begin{figure*}[h!]
\begin{center}
\includegraphics[angle=0,scale=0.52]{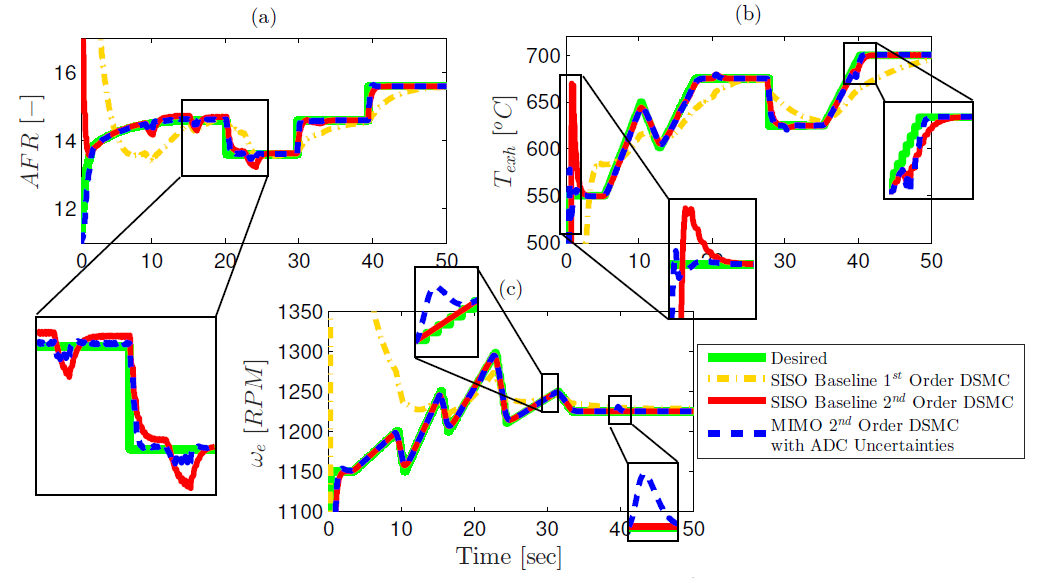} \vspace{-0.35cm}
\caption{\label{fig:All_MIMORobust_200ms_16bit}Results of desired trajectories tracking from SISO Baseline 1$^{st}$ and 2$^{nd}$ order DSMCs, and MIMO 2$^{nd}$ order DSMC with predicted ADC uncertainties ($T$=200~$ms$, quantization level=16-$bit$).} \vspace{-0.75cm}
\end{center}
\end{figure*}

Upon increasing the sampling rate from $20~ms$ to $80~ms$, and changing the ADC quantization level from 16-$bit$ to 10-$bit$, the first order DSMC performance degrades significantly. On the other side, the second order DSMC still presents accurate tracking results. By comparing the first and second order DSMC results at $T=80~ms$ and quantization level of 10-$bit$, it can be concluded that the proposed second order DSMC offers higher robustness against ADC uncertainties, and outperforms the first order controller by up to 90\% in terms of the mean tracking errors.
\begin{table} [t!]
\small
\begin{center}
\caption{Mean ($\bar{e}$) of Tracking Errors. Values Inside the Parentheses Show the Resulting Improvement from the Second Order DSMC Compared to the First Order DSMC.} \linespread{1.25}
\label{table:tracking_Results}\vspace{-0.15cm}
\renewcommand{\arraystretch}{0.85}
\begin{tabular}{lccccc}
        \hline\hline
\multicolumn{1}{c}{} & \multicolumn{2}{c}{$\bar{e}$~($T$=20$~ms$,~16-$bit$)}                                                        &  & \multicolumn{2}{c}{$\bar{e}$~($T$=80$~ms$,~10-$bit$)} \\
\cline{2-3} \cline{5-6}
\textbf{}& {$1^{st}$-Order} & {$2^{nd}$-Order} &  & {$1^{st}$-Order}  & {$2^{nd}$-Order} \\
\textbf{}& {DSMC} & {DSMC} &  & {DSMC}  & {DSMC} \\
\textbf{}   & {\textcolor{blue}{Reference}}& \textbf{}  &  & {\textcolor{blue}{Reference}} & \textbf{}  \\ \hline
AFR& 0.03  &  0.01   &    &  0.28   & 0.03  \\ \vspace{0.05cm}
[-]        &        & \textcolor{blue}{(-66.67\%)} &  &  & \textcolor{blue}{(-89.29\%)} \\ \hline \vspace{0.05cm} 
$T_{exh}$ & 0.2  & 0.1  &  & 4.0 & 0.4 \\ \vspace{0.05cm}
[$^o$C]                 &      & \textcolor{blue}{(-50\%)} &  &  & \textcolor{blue}{(-90.0\%)}\\ \hline \vspace{0.05cm} 
$N$  & 0.1  &   0.1   &  &   13.8    & 0.9\\ \vspace{0.05cm}
 [RPM]                &      & \textcolor{blue}{($\approx$0\%)} &  &  & \textcolor{blue}{(-93.5\%)} \\ 
\hline\hline 
\end{tabular} \vspace{-0.75cm}
\end{center}
\end{table}
\linespread{1} 

For the SI engine case study, $AFR$ controller is the most uncertainty-sensitive controller, in comparison with the engine speed and exhaust gas temperature controllers~\cite{Amini_SAE,AminiSAE2016}. Fig.~\ref{fig:All_MIMORobust_200ms_16bit} shows the tracking results of the first (SISO) and second (SISO and MIMO) order DSMCs under a relatively large sampling rate of 200$~ms$, which causes significant uncertainty at the controller I/O. As shown in Fig.~\ref{fig:All_MIMORobust_200ms_16bit}, the first order DSMC fails to track all the desired trajectories under these extreme ADC uncertainties, but the SISO second order DSMC shows acceptable tracking performances. 

Fig.~\ref{fig:All_MIMORobust_200ms_16bit} shows that among the three SISO second order controllers, the AFR controller is deviated more upon increasing the sampling rate. This deviation from the desired AFR trajectories is larger when there is a change in the desired engine speed trajectory (Fig.~\ref{fig:All_MIMORobust_200ms_16bit}). The link between engine speed and AFR controllers can be traced in the strong coupling between AFR and rotational dynamics via the intake air mass flow term ($\dot{m}_{ao}$):
\vspace{-0.45cm}
\begin{gather}
\label{eq:Engine_IEEE_4}
m_a(k+1)=m_a(k)+T(\dot{m}_{ai}(k)-\dot{m}_{ao}(k)).
\end{gather}
where, $\dot{m}_{ao}=AFR\dot{m}_f$. Moreover, Fig.~\ref{fig:All_MIMORobust_200ms_16bit} shows that the changes in desired engine speed profile have an effect on $T_{exh}$ controller performance. Similar to the AFR controller, the link between $T_{exh}$ and engine speed controllers is the exhaust gas temperature time constant ($\tau_e$) which is calculated with respect to the engine speed as $\tau_e=\frac{2\pi}{\omega_e}$.

The coupling within the engine dynamics can be represented in the MIMO DSMC design via $\boldsymbol{\beta}$ matrix. The diagonal elements of $\boldsymbol{\beta}$ are the same as the SISO controller, while the off-diagonal element represents the coupling between various sliding variables. According to Eq.~(\ref{eq:Engine_IEEE_4}) and $\tau_e$, $\boldsymbol{\beta}$ is defined to present the coupling between AFR, $T_{exh}$, and engine speed controllers. Additionally, the engine speed and air mass flow dynamics are inherently coupled because of the synthetic $m_{a,d}$ control input which is the input to the engine speed controller, and the reference trajectory for the air mass flow controller. Previously, the results in~\cite{stefanopoulou2000variable} showed that allowing $\dot{m}_{f,c}$ (used to regulate $AFR$) to depend upon the cam phasing (rotational dynamics) leads to smaller transients in $AFR$ tracking results. With a similar trend to~\cite{stefanopoulou2000variable}, here, Fig.~\ref{fig:All_MIMORobust_200ms_16bit} shows that by utilizing the MIMO controller with predicted ADC uncertainties, not only the AFR tracking error decreases (by 46\%), but also the effects of the engine speed trajectory variation on the AFR tracking become smaller. 

In a similar manner to $AFR$ controller, by utilizing the MIMO second order controller for $T_{exh}$, the desired $T_{exh}$ tracking performance improves slightly (by 11\%) compared to the SISO second order DSMC, and the large overshoot at the beginning of the simulation will be removed. On the other side, by looking into Fig.~\ref{fig:All_MIMORobust_200ms_16bit}, it can be seen that changing the engine speed controller to MIMO structure has almost no significant effect on the tracking performance. As highlighted in Fig.~\ref{fig:All_MIMORobust_200ms_16bit}, the MIMO controller results in small spikes, which occur when the desired AFR trajectory has a sudden change (Fig.~\ref{fig:All_MIMORobust_200ms_16bit}). \vspace{-0.5cm} 

\subsection{Handling Model Uncertainties} \label{sec:Engine_Control_adaptive}
In the next step, the performs of the proposed adaptation mechanism for the SISO/MIMO second order DSMC is investigated under up to 50\% multiplicative uncertainty within the engine model. In the following, the adaptation law for each of the engine controllers, along with the physical interpretation of the considered uncertain terms are discussed. 

\vspace{0.15cm}
$\bullet {\textbf{~~Exhaust~Gas~Temperature~Controller:}}$~According to Eq.~(\ref{eq:Engine_IEEE_2}), $T_{exh}$ dynamics is a strong function of the exhaust gas time constant ($\tau_{e}$). Thus, any error in estimating the $\tau_{e}$ directly affects the $T_{exh}$ dynamics and causes deviation from the nominal model. Error in estimating $\tau_{e}$ is represented by the multiplicative uncertainty term ($\alpha_{T_{exh}}$). The error in the modeled $T_{exh}$ dynamics is removed by using the following adaptation law with respect to Eq.~(\ref{eq:D2SMC_14}): \vspace{-0.45cm}
\begin{gather}
\label{eq:adaptive_Texh}
\hat{\alpha}_{T_{exh}}(k+1)=\hat{\alpha}_{T_{exh}}(k)+\frac{T(s_1(k))}{\tau_e \rho_{\alpha_1}}(600AFI-T_{exh}(k)).
\end{gather}

\vspace{0.05cm}
$\bullet { \textbf{~~Fuel~Flow~Rate~Controller:}}$~The fuel evaporation time constant $\tau_{f}$ plays an important role in the fuel flow rate dynamics. This means that any error in estimating $\tau_{f}$ can deviate the model under the test from the nominal model considerably. $\alpha_{\dot{m}_f}$ can represent any potential errors or variations in the estimated $\tau_{f}$. The adaptation law for $\alpha_{\dot{m}_f}$ becomes:  
\vspace{-0.45cm}
\begin{gather}
\label{eq:adaptive_mdotf}
\hat{\alpha}_{\dot{m}_f}(k+1)=\hat{\alpha}_{\dot{m}_f}(k)-\frac{T(s_2(k))}{\tau_f \rho_{\alpha_2}}\dot{m}_f(k),
\end{gather}
where, $\dot{m}_{f,d}$ is calculated according to desired AFR in Eq.~(\ref{eq:Engine_IEEE_4}).

\vspace{0.15cm}
$\bullet { \textbf{~~Engine~Speed~Controller:}}$~In the engine speed dynamics, the torque loss on the crankshaft ($T_{loss}$), is defined as $T_{loss}=0.4\omega_e+100$. The $T_{loss}$ is estimated by reading a torque map. Thus, the multiplicative uncertainty term $\alpha_{\omega_e}$ compensates for any error in reading the torque map. The adaptation law for the engine speed controller is:\vspace{-0.45cm}
\begin{gather}
\label{eq:adaptive_we}
\hat{\alpha}_{\omega_e}(k+1)=\hat{\alpha}_{\omega_e}(k)-\frac{T(s_3(k))}{J.\rho_{\alpha_3}}(0.4\omega_e(k)+100).
\end{gather}

\vspace{0.15cm} 
$\bullet { \textbf{~~Air~Mass~Flow~Controller:}}$~Air mass flow into the cylinder ($\dot{m}_{ao}$) is determined by $\dot{m}_{ao}=k_1\eta_{vol}m_a\omega_e$~\cite{Shaw}. $\eta_{vol}$ is the volumetric efficiency, and is calculated by using an empirical curve fit~\cite{Amini_CEP}. Thus, $\alpha_{m_a}$ can represent the uncertainty in $\dot{m}_{ao}$ that is extracted from $\eta_{vol}$ curve fit. $\hat{\alpha}_{m_a}$ is estimated by using the following adaptation law:  \vspace{-0.45cm}
\begin{gather}
\label{eq:adaptive_ma}
\hat{\alpha}_{m_a}(k+1)=\hat{\alpha}_{m_a}(k)-\frac{T(s_4(k))}{\rho_{\beta_4}}\dot{m}_{ao}.
\end{gather}

{As long as the adaptation algorithm parameter ($\rho_{\alpha}$) is tuned properly, there is no extra condition for convergence of the unknown parameters within the model. This has been examined and confirmed through extensive simulations for different trajectories and operation conditions, as shown in Sec.~\ref{sec:Results}} In Sec.~\ref{sec:Results}, the performance of both SISO and MIMO adaptive second order DSMCs will be evaluated in real-time by testing the controller software on a real ECU in a processor-in-the-loop (PIL) setup. \vspace{-0.4cm}

\section{Real-Time 2$^{nd}$ Order DSMC Verification} \label{sec:Results}
The designed adaptive SISO and MIMO second order DSMCs in Sec.~\ref{sec:Engine_Control_adaptive} are tested in a PIL setup (see~\cite{Amini_CEP} for the schematic and details of the PIL setup), to verify the performances of the adaptive DSMC in real-time. The PIL setup has two different processors, namely National Instrument (NI) PXI processor (NI PXIe-8135) and dSPACE MicroAutobox II (MABX). The model of the engine plant is built into the PXI processor. The output of the PXI processor is the controller feedback signal from the plant. Using embedded ADC units, the feedback signal is sampled and quantized at 80 $ms$ and 16-$bit$, respectively. On the other side, the adaptive second order DSMC logic along with the adaptation and uncertainty prediction mechanisms are implemented into the MABX, which is the main ECU. The output of the MABX is the control signal which is set to be updated at every 80 $ms$.
 
NI VeriStand$^{\textregistered}$~and dSPACE Control Desk$^{\textregistered}$~software on an interface desktop computer are used to configure the PIL setup, and conduct real-time tests, including desired trajectories tracking, unknown parameters estimation, and engine operation. Fig.~\ref{fig:alpha_adaptation_SiSO_80ms_16bit} shows the results of four unknown multiplicative parameters estimation from SISO and MIMO second order DSMCs with predicted ADC uncertainties. It can be seen that under up to 50\% uncertainty on each of the engine model's dynamics, by using the proposed adaptation mechanism, the unknown terms converge to their nominal values, ``1", in less than 4 $sec$. Moreover, the MIMO controller shows faster convergence for $\alpha_{T_{exh}}$ and $\alpha_{\dot{m}_f}$, while for $\alpha_{m_a}$ and $\alpha_{\omega_e}$, SISO and MIMO controllers show similar convergence behavior. 
\vspace{-0.7cm}
\begin{figure}[h!]
\begin{center}
\includegraphics[angle=0,width= \columnwidth]{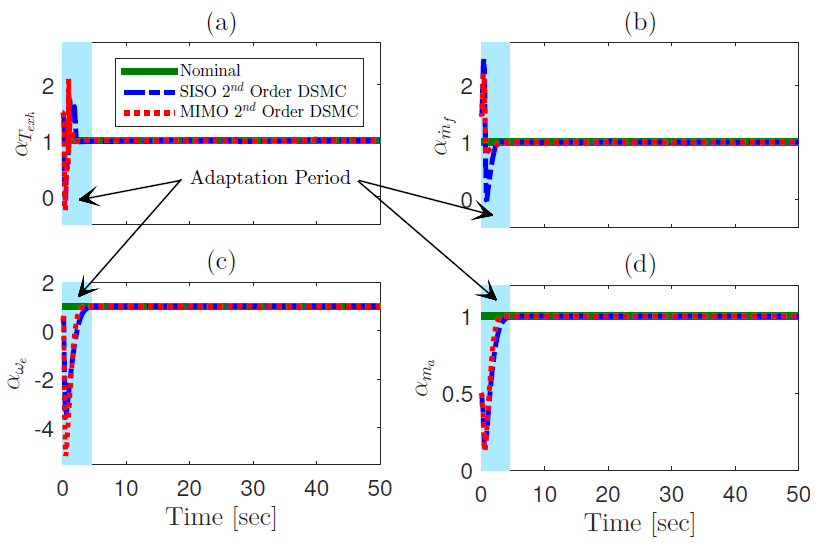} \vspace{-0.5cm}
\caption{\label{fig:alpha_adaptation_SiSO_80ms_16bit}Results of unknown multiplicative parameters convergences ($T$=80~$ms$, quantization level=16-$bit$).} \vspace{-0.85cm}
\end{center}
\end{figure}
\begin{figure*}[!h]
\begin{center}
\includegraphics[angle=0,scale=0.6]{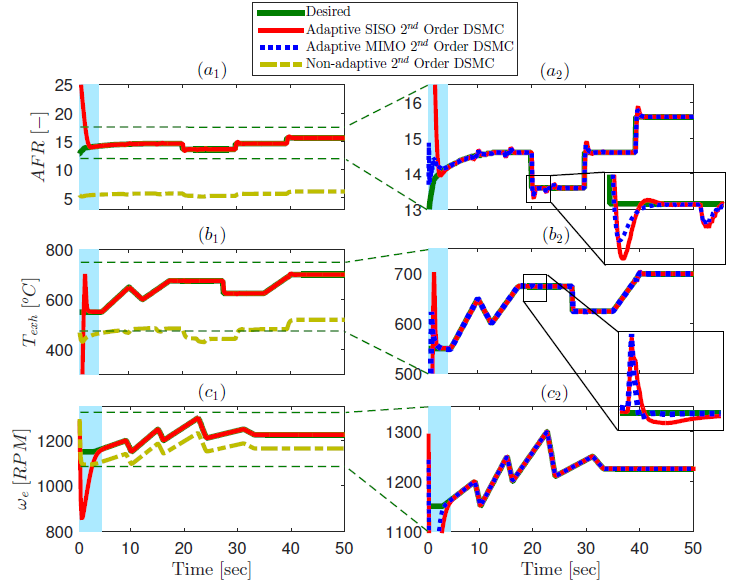} \vspace{-0.35cm}
\caption{\label{fig:adaptiveSISOMIMO_80ms_16bit}Results of engine control under model uncertainties: (a) air–fuel ratio, (b) exhaust gas temperature, and (c) engine speed ($T$=80~$ms$, quantization level=16-$bit$).} \vspace{-0.8cm}
\end{center}
\end{figure*}

Additionally, it can be observed that for all cases, despite the variation in the desired trajectories, after completion of the adaptation period, the uncertainties in the models are removed permanently. Finally, the PIL testing results show that the proposed adaptation mechanism is able to operate in real-time since it is computationally efficient. The adaptation time is significantly reduced if shorter sampling time is applied.

Fig.~\ref{fig:adaptiveSISOMIMO_80ms_16bit} shows the desired trajectory tracking results from the non-adaptive, and SISO/MIMO second order DSMCs. First of all, it can be seen that in the absence of the adaptation mechanism, due to the large uncertainties in the plant model (50\%), the non-adaptive controller fails to track the desired trajectories for all the cases. Upon activation of the adaptation mechanism, it can be seen that after completion of the adaptation period, both SISO and MIMO adaptive second order DSMCs with predicted ADC uncertainties provide accurate tracking performances. By comparing the non-adaptive and adaptive DSMCs it is revealed that the adaptation mechanism is able to remove the uncertainties in the model by more than 95\%, that consequently results in more than 90\% improvement in the controller tracking performance. 

Fig.~\ref{fig:adaptiveSISOMIMO_80ms_16bit} shows that the MIMO and SISO second order DSMCs have similar tracking behavior for the engine speed tracking. However, by using the MIMO structure, the tracking performance for $AFR$ and $T_{exh}$ controllers can be improved by 43\%, and 33\%, respectively. These improvements are more significant during the adaptation period, and also at those points where there are sudden changes in the desired engine speed profile. The latter observations can be explained based on the engine dynamics, in which the engine speed loop acts as a disturbance to $AFR$ and $T_{exh}$ (Eq.~(\ref{eq:Engine_IEEE_4})). By using the MIMO controller, $\dot{m}_{fc}$ and $\Delta$, which are used respectively to regulate $AFR$ and $T_{exh}$, can be configured to depend upon the engine speed dynamics. Linking the control input of the $AFR$ and $T_{exh}$ controllers to the rotational dynamics allows for better $AFR$ and $T_{exh}$ tracking performances during the engine speed transients.

{As another example, Fig.~\ref{fig:Rebuttal_fig} shows the closed-loop performance of the proposed second order adaptive DSMC in tracking non-smooth trajectories, under modeling and implementation imprecisions. The results in Fig.~\ref{fig:Rebuttal_fig} also show the evolution of the control signals, $\dot{m}_{ai}$, $\dot{m}_{fc}$, and $\Delta$. As can be seen, during the first few seconds, while the adaptation mechanism is estimation the unknown parameters, the control signals are rather aggressive, which results in deviations in the tracking. Once the unknown parameters estimation period is over in less than 5~$sec$, the trajectory tracking performance becomes nominal, and the control signals are smooth, except for those instants when a sudden change in the desired trajectory occurs.}   
\vspace{-0.75cm}
\begin{figure}[h!]
\begin{center}
\includegraphics[angle=0,width= \columnwidth]{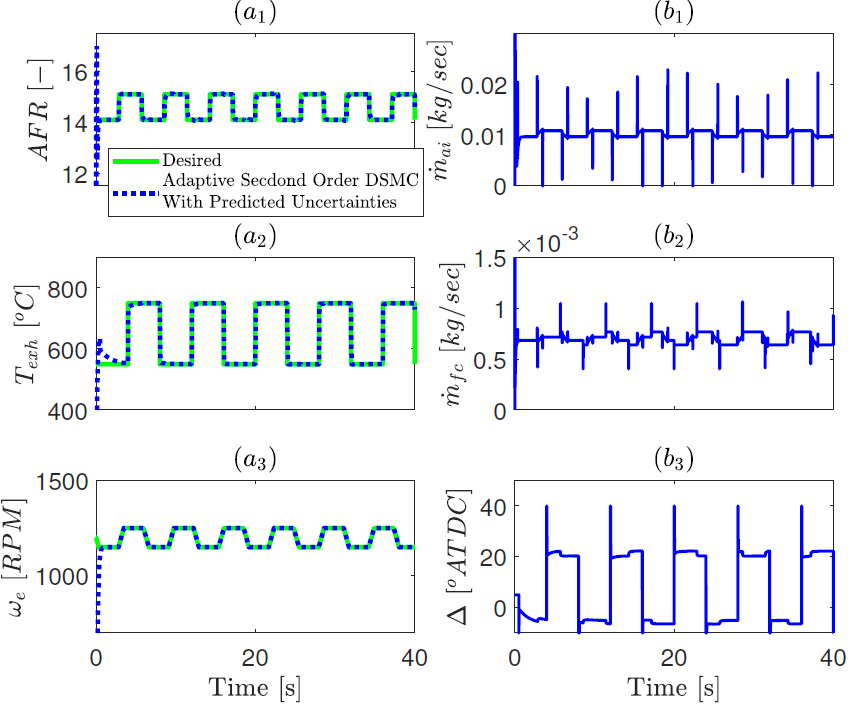} \vspace{-0.65cm}
\caption{\label{fig:Rebuttal_fig}{Performance of the adaptive 2$^{nd}$ order DSMC in tracking non-smooth trajectories: (a) Engine Performance, and (b) Control inputs. ($T$=20~$ms$, quantization level=16-$bit$)}.} \vspace{-1.0cm}
\end{center}
\end{figure}

\section{Summary and Conclusion}
\label{sec:conclusion}
A new formulation of an adaptive second order discrete sliding mode controller (DSMC) for multi-input multi-output (MIMO) nonlinear uncertain systems, along with an adaptation mechanism and a new switching control input, was presented in this paper. First, the adaptation law, for handling the uncertainties within the model, was driven based on a discrete Lyapunov stability theorem. Second, the behavior of the second order DSMC was studied on both reaching and sliding phases. In order to ensure the controller robustness against external analog-to-digital (ADC) imprecisions, a new switching control input was introduced, which contains the knowledge of ADC imprecisions via an online sampling and quantization uncertainties prediction and propagation mechanism. Third, the asymptotic stability of the proposed controller was guaranteed by invoking the new Invariance Principal for nonlinear discontinuous systems. 

The proposed controller was evaluated for a highly nonlinear combustion engine tracking control problem. The designed second order adaptive MIMO/SISO DSMC was tested in real-time on a real engine control unit (ECU) inside a processor-in-the-loop (PIL) setup. The simulation, and experimental real-time comparison results between the first and second order DSMCs revealed that:
\begin{enumerate}
\item The second order DSMC shows higher robustness against ADC uncertainties. When the sampling rate is changed from 20~$ms$ to 80~$ms$, and the quantization level changed from 16-$bit$ to 10-$bit$, the second order DSMC is able to improve the tracking errors by more than 90\% compared to the first order DSMC.
\item Inclusion of the physical coupling within the engine dynamics in the controller structure via the MIMO formulation allows for further improvement in $AFR$ and $T_{exh}$ controllers tracking performance (by up to 46\% and 11\%, respectively). However, the engine speed seems to lean strongly toward a decoupled or  minimally coupled structure.
\item In the presence of model uncertainties, it was shown that the proposed adaptation mechanism is able to remove the errors in the model permanently in less than 4~$sec$ of the engine operation time. Moreover, compared to the SISO DSMC, the MIMO controller shows faster unknown parameters convergence rates for the $AFR$ and $T_{exh}$ controllers. The MIMO adaptive second order DSMC is able to improve the $AFR$ and $T_{exh}$ tracking errors by 43\% and 33\%, respectively, compared to the SISO second order adaptive DSMC under modeling and ADC uncertainties. 
\end{enumerate}


\vspace{-0.75cm}
\begin{acknowledgment}
This material is based upon the work supported by the United States National Science Foundation under Grant No. 1434273. Prof. J. Karl Hedrick from University of California, Berkeley, Prof. Rush Robinett III from Michigan Technological University, and Dr. Ken Butts from Toyota Motor Engineering $\&$ Manufacturing North America are gratefully acknowledged for their technical comments during the course of this study.
\end{acknowledgment}
\vspace{-0.75cm}

\bibliographystyle{IEEEtran}

\bibliography{IEEEabrv,DSMC_UP_bib}

\end{document}